\documentclass[12pt]{amsart}
\title{la recensión de las leyes de Kepler en Inglaterra}
\author{jonathan Taborda}
\email{taborda50@gmail.com}
\dedicatory{A Johannes Kepler y Galileo Galilei, con motivo del cuadragésimo aniversario de la \textit{Astronomia Nova} y del \textit{Sidereus Nuncius}.}
\usepackage{amsmath,amssymb,amsthm,graphicx,graphics,cancel,pifont,pb-diagram,mathpazo}
\usepackage[activeacute,english,spanish]{babel}
\usepackage[T1]{fontenc}
\usepackage{lettrine,yfonts,calligra}
\newtheorem{teor}{Teorema}
\newtheorem{corol}{Corolario}

\begin{document}\maketitle
\renewcommand{\indexname}{Índice de Materias}
\begin{abstract}
El objetivo de la presente nota consiste en la discreta reconstrucción racional que se llevó acabo durante los años 1609-1630 y 1630-1666, i.e, el año de la publicación de la \textit{Astronomia Nova} y el año del fallecimiento del gran Astrónomo Alemán Johannes Kepler; y el posterior intervalo desde su deceso hasta el conocimiento de dichas leyes del movimiento planetario, por parte del Sabio Inglés Sir Isaac Newton.
\end{abstract}
\selectlanguage{english}
\begin{abstract}
The purpose of this note consists of discrete rational reconstruction which took place during the years 1609-1630 and 1630-1666, ie, the year of the publication of their \textit{Astronomia Nova}  and the year of death of the great German astronomer Johannes Kepler, and the subsequent range from desperate to knowledge of these laws of planetary motion, by the English scientist Sir Isaac Newton.
\end{abstract}
\selectlanguage{spanish}
\fontfamily{ppl}\selectfont
\lettrine{\calligra{R}}ecordemos que, este recorrido inicia en la antigua Grecia, con la
concepción Heliocéntrica de \textit{Aristarco de Samos}, la cual
pasa a través del triunfo de \textit{Ptolomeo} con su doctrina
Geocéntrica, y posterior declive y/o enterramiento en manos del
polaco \textit{N. Copernicus}.\\
Cuando apareció el \textit{De Revolutionibus} de Copernicus en 1543,
la disputa principal en Astronomía se daba entre los Averroístas,
quienes negaban la realidad para los epiciclos y excéntricas basadas
sobre argumentos de la física Aristotélica, y los Astrónomos
matemáticos, quienes soportaron la tradición teórica del libro de
texto y observaron los epiciclos y excéntricas como indispensables
para producir las posiciones de los cuerpos celestes. Los
Averroístas insistieron en que los Cielos estaban divididos en una
serie de órbitas concéntricas, todas centradas sobre la Tierra. Los
Astrónomos matemáticos continuaron una construcción empleando el
\textit{New Theories of the Planets} (Theoricae novae planetarum,
compuesto alrededor de 1460, primero publicada en 1472) de George
Peurbach (1423-1461).\footnotemark\footnotetext{Cf.
\textit{Theological Foundations of Kepler Astronomy.} Peter Barker
and Bernard R. Goldstein. The History of Science Society. 2001. pp.
90. ss.}\\
Peurbach empleó órbitas excéntricas, algunas de las cuales
acarreaban pequeñas esferas presentando la función de epiciclos.
Esta combinación estuvo acarreada por orbes internos y externos de
tamaño diferente, así que el sistema para los orbes de cada planeta
estaría centrado en superficies internas ó externas sobre la Tierra.
Los planetas estarían embebidos en los orbes pequeños
correspondiendo al epiciclo y transportados físicamente a través de
los cielos por movimientos combinados para el conjunto de orbes
completo.\\
La primera mayor publicación de Kepler fue su \textit{Mysterium
Cosmograpicum} (1596). El roll para la religión no es cancelado
sino que se indicó muy brevemente en el título del libro, que no
fue bien trasladado. ''Mysterium Cosmograpicum'' fue usualmente
interpretado como ''secret of the Universe''. Sólo que ''secret'' es
una traslación insípida de mysterium. El término puede significar
''mystery'' ó ''secret'', excepto que su significado central en la
antig\"{u}edad fue ''sacred mystery'', los secretos se enseñarían a
los iniciados cuando ellos entraran en la religión oculta. Así el
título debería ser mejor interpretado como '' The Sacred Mystery of
the Cosmos''.\\
Como es bien conocido, Kepler\index{Kepler, Johannes} introduce una construcción geométrica
basada sobre los cinco sólidos regulares platónicos para defender el
sistema Copernicano. El prefacio al lector inicia:
\begin{quote}
''Yo propongo, lector, demostrar en éste pequeño libro que el mayor
Dios y Gran Creador, en la creación para el movimiento del mundo, y
el ordenamiento de los cielos, referido a aquellos cinco sólidos
regulares, bien conocidos de Pitágoras a Platón en nuestro tiempo, y
que él naturalmente colocara en los cielos, otras proporciones, en
plan (ratio) para otros
movimientos''.\footnotemark\footnotetext{Cf. Ibid. pp. 99. n. 95.}
\end{quote}
\textit{A New Astronomy Based on Causes, or Celestial Physics}
(Astronomia Nova $AITIO\Lambda O\Gamma HTO\Sigma$. sev physica
coelestis) apareció en Praga en 1609. El libro inicia con una serie
de capítulos en que los sistemas de Ptolomeo, Brahe y Copernicus son
considerados como modelos posibles que pueden contar los datos de
Tycho para la posición de Marte con extremada exactitud. Este es una
investigación a posteriori para la Astronomía en el sentido que
prevaleció antes de Kepler-el objetivo es recubrir el fenómeno, no
dar una relación causal para el movimiento planetario conocido
realísticamente. En el capítulo 16 Kepler introduce un modelo que
emplea un ecuante con excentricidad no-bisectada, que el llama su
''vicarius hypothesis'', esto es, una hipótesis que es empleada
provisionalmente hasta el descubrimiento de la hipótesis
verdadera.\\
Por el final de la segunda mayor parte para la \textit{New
Astronomy} (finalizando con el Cap. 21) Kepler estuvo apelando al
ahora celebrado error de los 8 minutos de longitud para eliminar los
modelos de Brahe, Ptolomaee y Copernicus, autorizando únicamente su
hipótesis vicaria como una posible narración para la posición
angular de Marte.\footnotemark\footnotetext{Cf. Ibid. pp. 107. ss.}\\
También debemos subrayar, que de acuerdo a Barker y Goldstein, ''por
tanto las dos leyes reales presentadas en la \textit{New Astronomy}
no son la primera y la segunda Ley, como nosotros las conocemos en
la actualidad, excepto que son la Ley de distancia-velocidad y la
Ley de reciprocidad''.\\
La ley de reciprocidad es también llamada la ''regla seno-inverso'';
Kepler llama a esta ''libración''. La ''Ley de Áreas'' para círculos
es introducida en el Cap. 40 como una aproximación para la Ley de
distancia-velocidad, excepto que Kepler nunca le da el status de una
''Ley''. En el cap. 59, Kepler deriva la elipse de sus dos leyes, él
apela a la Ley de distancia-velocidad. La relación correcta entre la
Ley distancia-velocidad y la Ley de Áreas no fue establecida por
Kepler hasta el \textit{Epitome Astronomiae Copernicanae}, ó
\textit{Epitome of Copernican Astronomy} (Linz: Tampachius.
1618-1621), donde él indica que se tienen dos componentes para el
movimiento que conducen a la Elipse: uno es perpendicular para el
radio vector del Sol al Planeta, y el otro es recíproco a lo largo
del radio vector del Sol al planeta. Esto modifica su explicación
previa en la \textit{New Astronomy} y es equivalente a la Ley de
Áreas: ''Entonces en orden para la forma [la órbita elíptica] dos
elementos son mezclados conjuntamente, como ya se demostró: un
elemento será de la revolución al rededor del Sol por razón de una
virtud solar; el otro será de la libración en dirección del Sol por
razón de otra virtud solar distinta de la
primera''\footnotemark\footnotetext{Kepler, \textit{Epitome of
Copernican Astronomy}, KGW, Vol. 7. pp. 377. Cf. Barker and Goldstein. Op. cit. pp. 110. n.
70}.
Un rol fundamental, es también desempeñado por \textit{Tycho de
Brahe} (1546-1601) un noble Danés que además de beodo estaba
destinado ha dedicar su vida entera a la observación astronómica. Su
interés astronómico fue inspirado por el Eclipse solar de 1560. Durante 1580 y principios de 1590, Tycho derivó parámetros inapropiados para la teoría solar (esencialmente de Hipparchus), y trabajó sobre la producción para un catálogo extremadamente exacto de 777 estrellas. Este fue el primer catálogo disponible en Europa que fue independiente al de Ptolomy's, y representó un decrecimiento para un orden de magnitud en los errores comparados con el catálogo que Copernicus había incluido en el \textit{On the Revolutions}.\footnote{Mientras la influencia del sistema de Tycho sobrevivió brevemente en Europa, ésta tuvo una larga vida en China. Esta fue exportada por los Jesuitas, muchos de los cuales tenían muy buenos conocimientos en astronomía y para quienes el movimiento de la Tierra era un Anatema. El sistema de Tycho aparece en la Imperial Encyclopedia de 1726. La astronomía heliocéntrica para el Oriente realmente no encontró continuadores en China hasta inicios del siglo XVIII. Parte de la razón para este es el hecho de que el concenso principal para los chinos con los fenómenos astronómicos estuvo en la producción de calendarios, y los principios por los que fueron hechos empleaban aritmética en lugar de geometría. En adición, se tuvo una transformación muy pequeña de ideas entre China y el Oriente durante el siglo XVIII. Cfr. \textit{From Eudoxus to Einstein. A history of Mathematical Astronomy.} C. M. Linton. Cambridge Univ. Press. 2004. Cap. III. pp. 164. n. 24.}\\
Aquí entra en escena, nuestro Astrónomo/Astrólogo \textit{Johannes
Kepler}(1571-1630).\\
\textbf{Mysterium Cosmographicum} (1596).\\
Copérnico fue el primero capaz de construir un sistema astronómico compatible al de Ptolomeo, cuando no mejor, tanto desde el punto de las matemáticas como de la simple observación. Incluso entre los astrónomos más tradicionales, que no estaban dispuestos a considerar las consecuencias cosmológicas de una hipótesis de la astronomía matemática, la obra de Copérnico no solo fue respetada sino utilizada.\\
Independientemente de los grandes temas discutidos por los astrónomos mas destacados, la práctica de la astronomía y la astrología dependía en buena medida de las tablas astronómicas que, junto con determinadas reglas de uso o cánones, permitían determinar las posiciones pasadas o futuras de los planetas.\\
<<Aún así la obra de Copérnico estimuló la colaboración de modelos matemáticos geo-heliocentristas, como los de Erasmus Reinhold, Christopher Rothman- con su propia evolución hacia la cosmología copernicana- N. Reymers Ursus, Helisaeus Roeslin o Tycho Brahe, que situaban la Tierra en el centro del Universo, con el Sol girando a su alrededor, y los planetas al rededor del Sol. Tras diversas reivindicaciones y disputas por la prioridad, no siempre claras en todos los casos, este tipo de sistema pasó a conocerse simplemente como tychónico, reconociendo la autoría o, en todo caso, la autoridad de Tycho.>>\footnote{Cf. \textit{Talento y Poder.} Antonio Beltrán Marí. Historia de las relaciones entre Galileo y la Iglesia católica. $1^{a}$ Ed. Octubre de 2006. Editorial Laetoli. Cap. 2. pp. 85 y ss.}

En el sistema de \textit{Copérnico},\index{Copérnico} se tienen seis planetas girando
alrededor del Sol. Para el joven Kepler, tal sistema es una creación
perfecta de Dios, el hecho de que existan exactamente seis planetas
(éste no es el caso actual) debe tener una profunda razón.\\
Esta profunda razón fue uno de los misterios del Universo que el
joven Luterano decidió investigar. De acuerdo al mismo Kepler, una
revelación maravillosa le ocurrió el 19 de Julio de 1595, a saber,
la razón por la cual existen seis planetas, es por que en efecto,
existen cinco sólidos Platónicos (i,e. poliedros regulares) y cada
uno de ellos está en un lugar entre las seis "órbitas esféricas" tal
que éste poliedro está inscrito (ó circunscrito) en una de las seis
"órbitas esféricas". Así, tal estructura geométrica maravillosa no
únicamente se extiende exactamente a seis planetas, sino que ésta
también determina las razones entre el radio para las seis órbitas
esféricas.\\
Este es el origen para el primer libro de Kepler y él dedicó su vida
entera a estudiar el movimiento planetario con el objetivo de
verificar su ''conjetura salvaje".\\
En este texto, queda expuesta su teoría poliédrica y estructura del
universo; la cual refleja su Pitagorismo y Platonismo, pues aquí,
hace manifiesto su anexión a Proclo y revela su idolatría
Heliocéntrica, ya que esta gran obra de Astronomía-Astrología,
contiene el germen de la primera mitad de su otra obra cumbre, donde
queda probado geométricamente que ésta hipótesis de juventud, no es
tan descabellada, y su corroboración se da en el descubrimiento de
las armonías celestes.\\
''Como diría el propio Kepler en 1602, donde hay materia, hay
geometría; el mundo, geometría hecha materia, se encuentra dispuesto
según cantidades y relaciones geométricas. El mundo es la
realización de esa Idea, de ese arquetipo geométrico, en la materia:
la geometría, que es  Dios mismo, se manifiesta, se hace mundo, por
medio de la materia.\\
Por ello el mundo es esférico. Sólo así le fue posible a Dios
expresar la Trinidad en símbolo de ésta, con la <imagen de Dios
Uno-Trino en la superficie de la esfera; esto es, del Espíritu en la
regularidad de la ''relación'' entre el punto y la
circunferencia>."\footnotemark\footnotetext{Cf.Jorge M. Escobar. La mente de dios. Un estudio sobre la filosofía natural de johannes kepler. Cap. 3. Master's thesis,
Instituto de Filosofía. Universidad de Antioquia. 2003.}\\
Una segunda edición para el The Sacred Mystery of the Cosmos (con extensivas anotaciones) fue publicada en 1621, después de que Kepler había hecho sus descubrimientos sobre las tres leyes del movimiento planetario. Que los períodos orbitales dependían sobre las distancias al Sol dilucidó Kepler por un tiempo- él simplemente no tenía un dato lo suficientemente exacto para encontrar la ley correcta\footnote{Kepler en este tiempo consideraba que la relación entre el período orbital $T$ y la distancia $r$ era $T\propto r^2$, solo que esto no es cierto. En el \textit{Sacret Mystery} Cap. XX, Kepler propuso la relación $\left(T_2-T_1\right)\diagup T_1=2\left(r_2-r_1\right)\diagup r_1$, donde $r_1$ y $r_2$ son las distancias para los planetas sucesivos $(r_2>r_1)$ y $T_1$ y $T_2$ los respectivos períodos. Esta ley es equivalente a $\dfrac{T_2}{T_1}=2(r_2/r_1)-1$ que no es de la forma $T\propto r^\alpha$ para algún $\alpha$. Kepler posteriormente hipotetizó la ley $T\propto r^2$ en el Cap. 39. para la \textit{New Astronomy}. Cfr. \textit{From Eudoxus to Einstein}. Op. cit. pp. 173. n. 46.}-excepto que en retrospectiva conocemos que él se estaba formulando cuestiones más profundas.
En 1611, con su \textit{Strena de nive sexangula}, Kepler ampliaría
sus ideas sobre la realización del conjunto materia-arquetipo, pues
son Dios mismo, el mundo sólo empieza a existir con la materia; y
describirá con ella no sólo la creación del mundo en cuanto
totalidad astronómica y cosmológica sino también la generación de
los diversos seres que se pueden encontrar en él, tal como los copos
de nieve, los frutos de los árboles o las celdas de una colmena: en
todos los casos, se hallará involucrado un arquetipo que se realiza
en la materia.\footnotemark\footnotetext{Cf. La mente de dios. Op.cit. Cap. 3.}\\
En su último Opus magnum, \textit{Las Tabulae Rudolphinae} (1627),
Kepler asegura que la ciencia de los cielos se compone de dos
partes: la astronomía, que se ocupa de los movimientos de los
astros, cuyas leyes son inmutables, establecidas sobre los más altos
principios, y la astrología, que se ocupa, a partir de conjeturas,
de los efectos que tienen dichos movimientos en el dominio
sublunar.\\
En definitiva, si bien el \textit{Mysterium} atestigua los intereses
tempranos de Kepler en la astrología, no resulta útil como fuente al
momento de comprender el verdadero significado de la astrología
Kepleriana dentro de su obra y su pensamiento. Para ello debemos
referirnos a trabajos posteriores y, para empezar, al \textit{De
Fundamentis Astrologiae Certioribus} (1602), donde por primera vez
da una exposición sistemática de lo que piensa acerca de ella.\\
El \textit{De Fundamentis} Es un tratado astrológico y un calendario
para el año 1602, escrito poco después de la muerte de Tycho de
Brahe (24 de octubre de 1601).\footnotemark\footnotetext{Ibid.}\\
\textbf{Astronomia Nova}(1609)\footnotemark\\
\footnotetext{\textit{Astronomia Nova AITIO$\Lambda$ O$\Gamma$
HTO$\Sigma$, seu Physica Coelestis, tradita commentariis de Motibus
Stellae Martis, Ex observationibus G.V.Tychonis Brahe}, Prague,
1609, ahora disponible en dos excelentes ediciones por Max Caspar
(en Latín en \textit{Johannes Kepler: Gesammelte Werke}, Munich,
1973, iii; y en Alemán en \textit{Die Neue Astronomie}, Munich,
1929).} Kepler, de hecho, envía una copia del anterior libro a Tycho
de Brahe, \index{Brahe, Tycho} y análogamente, tal maestro de la astronomía deseaba
disuadirlo con respeto a que tal ''conjetura salvaje'' era una
fantasía de juventud. Sin embargo, él estuvo impresionado por la
inteligencia y originalidad de éste joven astrónomo. Para la década
de 1600, Tycho necesitaba el talento matemático del joven Kepler
para el ''conocimiento'' de sus observaciones astronómicas de toda
una vida, mientras que Kepler necesitaba el acceso a los datos
astronómicos de Tycho para verificar su gran ''misterio del
universo''.\footnotemark\footnotetext{Cf. \textit{The epic journey
from Kepler's laws to Newton's law of universal gravitation revited}
Hai-Chau Chang, Wu-Yi Hsiang. arXiv:0801.0308v1 [physics.hist-ph] 1
Jan 2008}\\
El sistema de Tycho fue defendido como una cosmología sólo tras la condena del copernicanismo. Sabido es que sus mayores defensores, los jesuitas, hicieron suyo el sistema de Tycho como consecuencia de su ciega obediencia a la Iglesia, no como resultado de la calidad de su ciencia. No hace falta decir que, en el caso de copernicanos como Kepler ó Galileo, el sistema de Tycho nunca fue tomado en serio como cosmología.\\
Kepler lo expresó muy gráficamente en una carta a Magini de Junio de 1601. En ella comenta la sencilla y <<simplísima>> forma en que el sistema copernicano explica los movimientos planetarios, recorriendo <<un círculo exactísimo en una revolución regularísima>>, exponiendo la relación entre períodos y las distancias de las órbitas:
\begin{quote}
\textit{<<no por fantasía sino por la verdad de las cosas. En realidad, Tycho también usa un ecuante en la Luna. Te darás cuenta fácilmente de hasta que punto la fuerza de Copérnico procede de esta armonía y simplicidad, esto es, de la perfección de los movimientos celestes. Pues por mas que Tycho imita y reproduce de cerca a Copérnico, retenida la Tierra inmóvil en el centro, no puede evitar que [la Tierra] refuerza en roscas de modo desigual y siempre diferente las vías por las cuales giran los planetas en el liquidísimo éter.>>}\footnotemark\footnotetext{Cf. Ibid. pp. 87.}
\end{quote}
Tycho Brahe esperaba que Kepler utilizara su precioso tesoro de observaciones para consolidar su sistema, pero finalmente redundó en beneficio del copernicanismo de Kepler. Por otro lado, cobraba fuerza la hipótesis del Sol como motor del movimiento planetario. <<Fue Kepler quien convirtió esta sugerencia en una tesis fundamental del copernicanismo y de la nueva filosofía celeste. Kepler convirtió al Sol en motor del movimiento planetario, primero en una formulación animista, y despues en clave mas mecanicista y ligada a una ley matemática.>>\footnotemark\footnotetext{Cf. Ibid. n. 4}
Esta estratégica unión, fructificó en el descubrimiento de la
primera y segunda ley del movimiento de Marte, a saber\\
\begin{enumerate}
\item [Primera ley:]\index{primera ley de Kepler} Marte se mueve describiendo una órbita elíptica
con el Sol en uno de sus focos.\\
\item [Segunda ley:]\index{segunda ley de Kepler} El radio vector trazado desde Sol a Marte barre
áreas iguales en tiempos iguales.
\end{enumerate}
El título completo en Latín, que es \textit{Astronomia Nova...}, se trasladó como \textit{Nueva Astronomía Basada sobre causas, ó Física celeste, tratada por el significado de los comentarios sobre los movimientos para las estrellas de Marte a partir de las observaciones de Tycho Brahe}. La primera traslación Inglesa fue publicada recientemente en 1992. Algunos pasajes están trasladados en Koyré (1973) en el que el trabajo es discutido con algún detalle. Otro texto en el que se da una narración para los contenidos de la \textit{New Astronomy}- particularmente los aspectos de la astronomía de Kepler- es Stephenson (1987) y otros aspectos técnicos para el trabajo son descritos entre otros en Aiton (1969) y Whiteside (1974).\footnote{Cf. Linton. Op. Cit. n. 56.}\\
El estilo de la \textit{New Astronomy} fue totalmente diferente a el Ptolomy ó Copernicus. No únicamente fue el estilo novedoso de Kepler, la Astronomía que él presentó era realmente nueva, con un énfasis completamente diferente al de sus predecesores. En lugar de intentar reproducir los movimientos celestes como apropiados posiblemente con ciertas construcciones geométricas especiales, Kepler intentó conocer la estructura básica del Universo de Dios.\\
Kepler entonces describió su propia teoría de la <<gravitación>>, que contiene los gérmenes de lo que sería, en manos de Newton, la teoría de gravitación universal; él pensaba la gravitación como una atracción mutua entre cuerpos similares para la fuerza magnética. El método de Kepler para obtener con exactitud la posición de la Tierra fue ingenioso. Él determinó la posición para la Tierra de las observaciones de Marte, proporcionando el lugar de Marte con exactitud. Las longitudes de Marte fueron bien establecidas, como es bien observado, exceptuando sus distancias. Kepler va alrededor de estas por el empleo de observaciones para Marte separadas por su período zodiacal (687 dias) así que, sinembargo para Marte fue el del Sol, y este era el mismo tiempo. Una vez obtenido lo anterior él lo administró para mostrar que la órbita de la Tierra era mejor representada por el mecanismo del ecuante-excentricidad Ptolemaico, similarmente para otros planetas.\\
Kepler apreció la equivalencia entre las teorías planetarias Copernicanas y Ptolemaicas, y él conocía que, en el esquema de Ptolomeo, el epiciclo representaba la órbita de la Tierra al rededor del Sol.\\
Luego de la dura batalla librada por Kepler, en el establecimiento de la trayectoria para Marte, y después de modificar su trayectoria oval, conocemos que de su primera ley del movimiento planetario la palabra <<foco>> no fue mencionada en el corpus central del \textit{New Astronomy}, y fue únicamente hasta su \textit{Epitome of Copernican Astronomy} que Kepler enfatizó este aspecto para las órbitas planetarias. El mecanismo físico de Kepler por el que los planetas fueron forzados a moverse de tal forma fue incorrecto, exceptuando la idea de que el estudio para los movimientos celestes estaría basada sobre causas físicas influenciando todos aquellos que siguieron sus pasos.\\
La mayoría de los contemporáneos de Kepler observaron no convincente la condición de tornar la astronomía en una ciencia física. Incluso Michael Mastlin tuvo sus propias dudas. En una carta Kepler, de 1616 escribió:<<Concerniente al movimiento para la Luna tu escribes que tuvisteis que trazar todas las desigualdades para las causas físicas; Yo no estoy completamente enterado de esto. Yo pienso que en lugar de tal uno debería dejar las causas físicas fuera de contexto, y explicar las materias astronómicas únicamente de acuerdo al método astronómico utilizado en tal, no físico, causas e hipótesis.>>\footnote{Cf. Ibid. pp. 190. n. 79.} Las primeras leyes del movimiento planetario fueron realmente revolucionarias. Con ellas él obtuvo la premisa fundamental de toda la astronomía matemática que le precedió, i.e., que los movimientos celestes fueron construidos sobre movimientos circulares uniformes. El hecho de que Kepler descubriera esas leyes es el acontecimiento más importante de todos por el acto de haberlo realizado no solo con los datos observacionales de Tycho solamente, y sin contar con la ayuda de la dinámica Newtoniana. Crucial para éste descubrimiento fue su insistencia de que el movimiento para los planetas fue un proceso físico, con el Sol jugando un rol predominante. Aunque la física de Kepler fue errónea, éste no dejó de formular algunas cuestiones importantes y responderlas correctamente.\\

\textbf{Harmonices Mundi} (1619)\\
En 1613, como es conocido, Kepler estuvo interesado en el problema de determinar los volúmenes para barriles de vino en tamaños diferentes. Él tomó las ideas que había desarrollado cuando determinó las áreas para regiones curvadas como parte de su guerra con Marte, este desembocó en la publicación de su \textit{New Solid Geometry of Wine Barrels} <<\textit{Nova steriometria doliorum vinariorum}>> en 1615. El método de Kepler involucró la división de volúmenes en el interior de muchas regiones infinitesimales y representó un primitivo cálculo integral. Este fue expandido sistemáticamente por Bonaventura Cavalieri. Kepler también planeó producir un libro de texto que combinara la teoría heliocéntrica de Copernicus con sus propios descubrimientos en un folio verdaderamente aceptable.\\
La primera parte de este trabajo, el \textit{Epitome of Copernican Astronomy} fue impreso en 1617, excepto que el libro no estuvo completo hasta 1621.\\
El \textit{Harmony of the World} es un libro fascinante, con discusiones matemáticas tales como la constructibilidad (usando regla y compás únicamente) para polígonos regulares, teselaciones para el plano y poliedros semi-regulares ó sólidos Arquimedianos (Kepler da la primera prueba conocida de que se tienen exactamente 30 de ellos.)\footnote{Un sólido Arquimediano es un poliedro convexo en el que todas sus caras son polígonos regulares (para a lo sumo dos tipos diferentes) y los vértices para tales todos idénticos. Cf. Ibid.} Kepler también discutió los méritos relativos para los métodos geométricos largamente establecidos y los nuevos desarrollos en técnicas algebraicas para la solución de varios problemas. V.g., Kepler creía que los polígonos regulares que no eran construidos empleando regla y compás eran <<no conocidos>>, y no jugaron un rol en la construcción del Universo como designio de Dios.\footnote{En tiempos de Kepler, los únicos polígonos regulares construibles conocidos eran aquellos que tenían 3, 4 y 5 lados y aquellos derivados de ellos, i.e., aquellos con el número de lados $2^n3^m5^i$ con $n\in\mathbb{Z}^+$ y $m, i=0\,\text{ó}\, 1$. En efecto, Kepler incluye una prueba en el libro I del \textit{Harmony} en el que los polígonos regulares con un número primo de lados superior a 5 no son construibles, solo que su prueba es incorrecta. En 1796, con solo 19 años de edad, Gauss probó que un polígono regular con un número primo de lados $p$ es construible sii, $p-1$ es una potencia de 2. El resultado de Gauss muestra que es posible construir polígonos regulares con 3, 5, 17, 257 y 65.537 lados. No se conoce si existen algunos otros primos p. Cf. Ibid. n. 88.}

Llevó a Kepler otra década de duro trabajo verificar que el mismo
principio de la primera y segunda ley también son fijos para los
otros seis planetas, y además, el descubrió la siguiente tercera
ley, a saber:\\
\begin{enumerate}
\item [Tercera ley:]\index{tercera ley de Kepler} La razón entre el cubo del eje mayor y el
cuadrado para el período, i.e. $\dfrac{(2a)^3}{T^2}$, es la misma
para todos los seis planetas.\footnotemark\footnotetext{Cf. arXiv:
0801.0308v1 [physics.hist-ph] 1 Jan 2008. pág. 6. ss.}
\end{enumerate}
Pero dejemos que sea el mismo Kepler, quien nos narre su asombro:
\begin{quote}
\fontfamily{pzc}\selectfont No han pasado ni dieciocho meses desde
que vi el primer rayo de luz, ni tres meses desde que amaneció, y
muy pocos días desde que el Sol, en todo su esplendor, lo más
admirable que se puede ver, brilló repentinamente ante mí. Nada me
detiene; no me voy a culpar por mi furia sagrada; triunfaré sobre la
humanidad cuando confiese honestamente que he robado los vasos de
oro de los egipcios para construirle un tabernáculo a mi Dios lejos
de los confines de Egipto. Si me perdonáis, me alegraré; si os
ponéis furiosos conmigo, podré soportarlo; la suerte está echada, el
libro está escrito para que se lea ahora o en el futuro. No me
importa quién lo lea; puede esperar un siglo hasta que surja un
lector, dado que Dios ha esperado seis mil años para que alguien
observara su obra.\footnotemark\footnotetext{Cf. \textit{De aquí al
infinito. Las matemáticas de hoy.} Ian Stewart. Drakontos. Crítica.
Traducción castellana de Mercedes García Garmilla. 1998. Cap. 14.}
\end{quote}
Aunque las ideas de Kepler revelan de forma inequívoca la influencia
de Paracelso y de sus discípulos, el antagonismo entre su método de
aproximación científico y la actitud mágico-simbólica de la alquimia
era, pese a todo, tan fuerte que Flud, famoso alquimista de su
tiempo y miembro de los rosacruces, mantuvo una violenta polémica a
propósito de la obra maestra de Kepler, \textit{Harmonices
mundi}.\footnotemark\footnotetext{Cf. \textit{Escritos sobre física
y filosofía.} Wolfgang Pauli. Debate. pensamiento. Ed. por Charles
P. Enz y Karl von Meyenn. Versión castellana de Mercedes García y
Rodolfo Hernández. 1996. Ensayo 21. Pág. 277, ss.}\\

Ahora, trataremos de narrar la interpretación dinámica de las
mencionadas Leyes del movimiento planetario de Kepler, desde la
metodología y/o estilo Newtoniano manifestado perogrullamente en el
\textit{Principia}. A grandes rasgos, los mayores resultados para el
análisis matemático de Newton en su magna obra pueden ser
establecidos como los siguientes teoremas, a
saber\footnotemark\footnotetext{Cf. Chang-Hisiang. Ibid. pág. 2. ss.}:\\
\begin{teor}
La segunda ley de Kepler es fija (i.e,
$\dfrac{dA}{dt}=\text{constante}$) sii la aceleración (res. la
fuerza) está en dirección al centro (i.e.
centripetal).\hspace{0.5cm}\qedsymbol
\end{teor}
\begin{teor}
La primera y segunda ley de Kepler implican que el vector
aceleración \textbf{a} está apuntando en dirección del Sol y con su
magnitud igual a\\$$\dfrac{\pi^2}{2}\dfrac{(2a)^3}{T^2}\dfrac{1}{r^2}$$\\
donde\\$$\textbf{a}=\dfrac{\pi^2}{2}\dfrac{(2a)^3}{T^2}\dfrac{1}{r^2}\left(\begin{array}{c}-\cos\theta\\-\sin\theta\end{array}\right)$$\hspace{0.5cm}\qedsymbol
\end{teor}
\begin{teor}[\textbf{De unicidad y recíproco para el Thm. 2} ]
Supóngase que el vector aceleración es centrípeto y con su magnitud
inversamente proporcional al cuadrado de la distancia, donde\\
$$\textbf{a}=\dfrac{K}{r^2}\left(\begin{array}{c}-\cos\theta\\-\sin\theta\end{array}\right)$$\\
Entonces, el movimiento satisface la segunda ley de Kepler y su
órbita es una sección cónica.\hspace{0.5cm}\qedsymbol
\end{teor}
\begin{teor}
La fuerza de gravedad para una forma esférica con densidad uniforme
(preservando el área) sobre una partícula externa P es igual
a\\$$G\dfrac{Mm}{{\overline{OP}}^2}$$\\ donde M (resp. m) es la masa
total para la forma esférica (resp. la masa para la partícula en P)
y $\overline{OP}$ es la distancia entre el centro O y
P.\hspace{0.5cm}\qedsymbol
\end{teor}
\begin{corol}
Sean $\Sigma_1$ y $\Sigma_2$ una par de cuerpos esféricos con
densidad radial uniforme (i.e. cada uno de ellos puede ser
descompuesto en el interior de la unión para las formas esféricas
del Thm. 4) Entonces la magnitud para la fuerza de gravitación
(total) entre ellas es igual a\\$$G\dfrac{M_1M_2}{{\overline{O_1O_2}}^2}$$\\
donde $M_1$ (resp. $M_2$) son las masas total para $\Sigma_1$ (resp.
$\Sigma_2$) y $O_1$ (resp. $O_2$) son sus
centros.\hspace{0.5cm}\qedsymbol
\end{corol}
Los anteriores teoremas, son una interpretación moderna de los
presentados por el Sabio Inglés, en el paper de Chang y Hsiang.\\
Sabido es que, el astrónomo Alemán, fue uno de los precursores en la
introducción del concepto de ''inercia'' en el discurso de la
filosofía natural; y que tal concepto, transformado en las manos de
Newton, \index{Newton, Isaac} fue fundamental en la batalla entre las concepciones
mecanicistas Cartesianas relativas al movimiento de los cuerpos en
medios resistentes (libro II del \textit{Principia}),\index{Principia} y cuyas
tres leyes del movimiento planetario desde la heurística Newtoniana,
constituyen el declive y/o ocaso absoluto con respecto a la teoría
de los vórtices Cartesianos.\\
Ahora bien, regresando a Kepler, el interesante artículo del Prof. Victor E Thoren, nos ayudará a
clarificar en gran medida la historia de la segunda Ley de Kepler en
Inglaterra. En efecto, el Paper del profesor Thoren, es un análisis
pormenorizado de otros dos Paper's; Russell:\textit{Kepler laws of
planetary motion, 1609-1666.} BJHS. (1964). pp. 1-24, y Witheside \index{Witheside, Derek Tom};
\textit{Newton early trought motion:
a fresh look}. BJHS. (1964). pp. 117-37.\\
Aunque concebidos independientemente y con objetivos diametralmente
opuestos, el profesor Thoren concluye que dichas investigaciones
coinciden únicamente en una área substancial, i.e, las citaciones
para la segunda Ley de Kepler por parte de los astrónomos Ingleses
datan entre 1650 y 1670.\\
Es bien conocido desde el punto de vista histórico-epistemológico,
el considerar a Kepler y a Newton como los dos grandes genios del
siglo XVII. Mientras los voluminosos escritos Keplerianos pudieron
tener una mayor circulación, esto es debido claramente a su peculiar
aprovechamiento físico para la Astronomía que fue recibido con menor
entusiasmo del que gozó el peculiar aprovechamiento matemático para
la física realizado por Newton.\\
Muchos de los Astrónomos Ingleses más respetados tales como Jeremiah
Horrox tenían  en gran estima los trabajos Keplerianos, en especial
el trabajo iniciado conjuntamente por él y el gran
observador Tycho Brahe; \textit{Rudolphine tables}. \index{Rudolphine tables}\\
En Vincent Wing, un gran pionero de la Astronomía Isabelina,
encontramos varias referencias a la segunda Ley de Kepler,
específicamente en su: \textit{Harmonicon Coeleste} (1651).\\
Por tanto, en Jeremy Shakerley en su \textit{Anatomy of Urania
Practica} de 1649, manifiesta también una gran evidencia de esto.
Escrito en un estilo animadvertido para la época, éste es una
crítica para Wing. Lo que es significativo en Shakerley es la
habilidad para citar el capítulo y el verso de Kepler sobre algún
tópico. Sí él no estaba extremadamente familiarizado con Kepler, se
tuvo que haber visto obligado a dar la impresión de que lo estaba.\\
Wing también, por ésta fecha, se refería únicamente a las
\textit{Tablas Rudolfinas} en sus efemérides. Excepto que cuando
publicó su propio texto en 1651, él empleó un ecuante de foco-vacío
modificado, que él racionalizó como sigue:
\begin{quote}
\textit{''Ahora para hallar el lugar de un planeta en su Elíptica
Kepler en su Epit. Astro. Copernic. enseñó como puede ser obtenida,
excepto que Bullialdus (para hacer la operación más sencilla)
muestra la presentación para el mismo por un Epiciclo, cuyo
movimiento es doble para el movimiento de un planeta en su órbita, y
así por la solución para las líneas rectas triangulares, éstas
pueden ser encontradas con más facilidad, en la dirección que ellas
toman, como muy racionalmente''...}\footnotemark\footnotetext{
Vincent Wing, \textit{Harmonicon Coeleste} (London, 1651). pp. 44.
Cf. Thoren. Op. cit. pp. 246. n. 13.}
\end{quote}
Shakerley aparentemente siente el mismo peso para sus lectores.\\
En sus \textit{Tabulae britannicae;...calculated...from the hypothesis
of Bullialdus,and the observations of Mr. Horrox of 1653} , fue más
condescendiente con sus lectores en la mención para una elipse,
refiriéndose en todas partes a las órbitas como ''excéntricas''.\\
Samuel Foster, profesor de Astronomía en el Gresham College, fue más
cándido:
\begin{quote}
\textit{El camino que yo tomo es (en general) el de agregarme a
Copernicus...y en particular, para con Kepler en sus Rudolphine
Tables. Únicamente se tiene una diferencia: Kepler hace que las
órbitas para los planetas sean elipses, que es la mejor dirección;
ó aquí hace que sean círculos perfectos, que es el camino más
sencillo.}\footnotemark\footnotetext{\textit{Miscellanies of Mr.
Samuel Foster} (London, 1659). Despista el dato de publicación,
Foster falleció en 1652. Cf. Ibid. n. 15.}
\end{quote}
En Oxford durante el siglo XVII apareció la primera manifestación
seria de astronomía académica en Inglaterra.\\
Centralizado sobre Wilkins, ésta actividad fue continuada por Seth
Ward y John Wallis, y soportada por figuras como Christopher Wren,
Lawrence Rocke o Paul Neile.\footnotemark\footnotetext{Todos excepto
Ward pertenecieron al pequeño grupo de miembros fundadores de la
Royal Society. Es posible que Hooke fuese un miembro de dicho grupo,
desde su estancia en Oxford en 1653 y a lo largo de su vida
manifestó serios intereses en la Astronomía. Su celebrada carta a
Newton del 6 de Enero de 1679/80 muestra que él estaba lo
suficientemente familiarizado con Kepler para citar su formulación
alternativa para la segunda Ley (a saber, que la velocidad para un
planeta en su órbita es inversamente proporcional de su distancia al
Sol.) Cf. Ibid. n. 17.}\\
En 1669 Cassini elaboró su teoría sobre excentricidades en su
''\textit{Methodus investigandi apogeae, excentricitates et
anomalias planetarum, breviter exposita et demonstrata}'', que
excepto por su tiempo y marco geométrico sonoro no fue suficiente:
la precisión astronómica también fue demandada. Sería Mercator quien
daría el paso hacia adelante para argüir que él tomó la estabilidad
del hecho de que el simple ecuante de foco-vacío era inapropiado en
la Astronomía de posición. Si nosotros creemos a Mercator, Kepler
tuvo que probar con claridad que la teoría del ecuante simple era
insostenible, algo que nunca hizo posteriormente por
entretenimiento.\\
El primer astrónomo Royal John Flamsteed quien mantiene una intensa
colaboración con Newton, tuvo una reacción inicial con Kepler
bastante desafortunada. Él rápidamente despachó las especulaciones
físicas como ''nociones infundadas'', y en su primer paper
publicado, sobre efemérides ''para los fenómenos más notables del
año 1670'', emplea a Streete y Wing para su trabajo.\\
Con Flamsteed instalado como Astrónomo oficial de la nación, Halley
fue el Astrónomo Inglés más importante durante el medio siglo
siguiente; entre ambos fundaron el aparato técnico que Newton
emplearía como contenido astronómico en el \textit{Principia}.\\
Así, mientras Newton estuvo hablando de sus propias experiencias ò
consultas, fue enteramente apropiado que él tuvo que
etiquetar la segunda Ley de Kepler como ''\textit{propositio
Astronomis notissima}''.\\
Como ya es conocido, la segunda Ley fue formulada originalmente en
1609, en dos formas diferentes: la forma correcta es llamada Ley de
áreas, y se tiene una forma alternativa que puede establecerse como
la velocidad para un planeta varía inversamente como su distancia al
Sol. Esta es denominada la Ley de distancia inversa. En un tiempo
Kepler las observó a ellas como matemáticamente equivalentes, pero
en efecto ellas no lo son; la Ley de distancia inversa es una buena
aproximación para elipses de excentricidad pequeña, excepto que no
son exactas. Para 1621, sin embargo, cuando la última parte para el
\textit{Epitome} fue publicada, él tuvo que probar que las dos leyes
no eran idénticas y que la Ley de áreas era correcta.\\
Cuando Kepler publicó su teoría en 1609 la rotación solar no había
sido aun observada. Él entonces postuló sobre argumentos puramente
teóricos un período de menos de tres meses- éste es el tiempo
empleado por Mercurio-, el acercamiento del planeta al Sol, para
completar su órbita. Cuando, pocos años después, las observaciones
de Galileo sobre las manchas solares mostraron que el Sol hace su
rotación en la dirección requerida, con un período de 28 dias,
Kepler naturalmente observó esto como una fuerte confirmación para
su teoría. Analizando el período perteneciente a 1609-1630,
encontramos que la \textit{Astronomia Nova} llamó poco la atención
cuando ésta fue publicada por vez primera. Este es un libro de
difícil lectura. Éste es difuso, y mucho de éste es simplemente un
recuento de los ensayos tempranos de Kepler para resolver el
problema de la órbita de Marte. No es hasta la pp. 284 (de las 337
que lo componen) que las dos primeras leyes aparecen. Kepler y sus
contemporáneos no estaban completamente familiarizados con las
propiedades de las elipses; su aprovechamiento matemático para ellas
es tosco y no sistemático, mientras que sus lectores menos
cualificados que él conocían las propiedades y las aplicaban para
computaciones astronómicas.\\
Para 1630, se han encontrado pocas referencias para las ideas de
Kepler.\\
Uno de los primeros lectores de la \textit{Astronomia Nova} fue el
Astrónomo y Matemático, Thomas Hariot, quien recibió una copia
prontamente después de la publicación y recomendó ésta a otro
Matemático y pupilo suyo, William Lower. Se conoce la reacción de
Lower por una carta que él le escribió a Hariot en Febrero de
1610.\footnotemark\footnotetext{Publicado en\textit{Thomas Hariot:
the Mathematician, the Philosopher, and the Scholar}, by Henry
Stevens;1900. pp. 120-124. Cf. \textit{Kepler Laws of planetary
motions: 1609-1666}. J.L.Russell. BJHS. Vol 2. $N^\circ 5$ (1964).
n. 8.}\\
Él claramente encontró el trabajo casi intolerablemente difícil,
excepto que al mismo tiempo intensivamente estimulante. De su
lectura él aceptó muchas de las ideas de Kepler, incluyendo las
órbitas elípticas, excepto que él posteriormente necesitó la ayuda
de Hariot. Parece que de la carta de Lower's, Hariot por sí mismo
aceptó las ideas de Kepler, a lo sumo en sustancia, excepto que él
no publicó ninguna cosa sobre el tema.\\
En 1612 el sabio Italiano, Federico Cesi, un amigo y patrón de
Galileo y joven miembro de la Lyncean Academy, en una carta escrita
para Galileo, mencionó la teoría planetaria Kepleriana de las
elipses con aprobación. Es importante mostrar que Galileo estuvo al
corriente de dicha teoría, puesto que él nunca mencionó esto en sus
escritos y ciertamente él no aceptó esto. El soporte más fuerte fue
en 1615 cuando Giovanni Magini, Profesor de Matemáticas en Bologna,
publicó su \textit{Supplementum Ephemeridum} en que él empleó las
leyes de Kepler en la computación de efemérides para Marte. Sin
embargo, aparte de un conocimiento general de como él aplicó la
teoría de Kepler, él no da detalles de como era la teoría.\\
La publicación para el \textit{Harmonices Mundi} \index{Harmonices Mundi} en 1619 esparció el
conocimiento de las ideas del autor, aun reducido número de
astrónomos. Éste fue aparentemente el menos leído de sus mayores
trabajos y se tienen pocas referencias para éste en la literatura
astronómica subsecuente. Su próximo trabajo, \textit{Epitome
Astronomiae Copernicanae} \index{Epitome
Astronomiae Copernicanae}, fue una poderosa defensa para el Sistema
Copernicano, en el curso que el da una narración completa para
ambas, su teoría física y sus tres leyes. Este fue publicado en tres
partes en 1618, 1620 y 1621, y eventualmente sería muy influenciable
y ampliamente leído, excepto que éste lo hizo lentamente famoso.\\
La teoría Kepleriana para Marte fue brevemente notificada por dos
escritores en 1622: Longomontanus en \textit{Astronomia Danica} y
Nathanael Carpenter en la segunda ed. para su \textit{Philosophia
Libera}. Ambos rechazaron las elipses; ellos aquí estaban mal
dispuestos para abandonar el principio de que los movimientos
planetarios serían reducidos a círculos perfectos.\\
Otras críticas se repitieron en la segunda ed. de la
\textit{Astronomia Danica} (1640) y la tercera ed. para la
\textit{Philosophia Libera} (1635). Kepler tuvo, sin embargo, a lo
sumo un discípulo durante los inicios de 1620: Philip M\"{u}ller,
profesor de Matemáticas en la Universidad de Leipzig. No se observa
que M\"{u}ller haya publicado algo sobre el tema, excepto su
aceptación para sus ideas que es mostrada en sus cartas para Kepler
y en su
correspondencia con Peter Cr\"{u}ger.\\
Se tiene alguna evidencia de que Willebrord Snel (1591-1626) también
aceptó las elipses. El astrónomo Alemán, Ambrosius Rhodius, en un
prefacio de Michael Havemann´s Astraea (1624), comentó la
\textit{Astronomia Nova} para su defensa del sistema Copernicano.\\
Él dispuso claramente bien de las ideas Keplerianas.\\
El punto de retorno en la fortuna Kepleriana lo constituyó la
publicación de su último gran trabajo, \textit{Rudolphine Tables},
en 1627. Este fue un evento por el que el mundo científico esperó
largo tiempo. Ellas están basadas sobre datos observacionales de
primera clase acumulados por Tycho Brahe en la última parte del
siglo XVII, y en ellas, por vez primera, las leyes fueron colocadas
en conjunto. A partir de ahí, los astrónomos compararían las
predicciones para las tablas con las posiciones actualmente
observadas para el Sol, la Luna y los planetas, y compararían
entonces los resultados con aquellas teorías astronómicas rivales.\\
La impresión hecha por esas tablas sobre Peter Cr\"{u}ger es
vívidamente conocida por extractos con Philip M\"{u}ller, publicada
por von Dyck y Caspar en 1927.\\
Únicamente la parte de Cr\"{u}ger está disponible, excepto que
claramente M\"{u}ller fue un soporte de Kepler y, de 1620 en
adelante, fue Cr\"{u}ger quién estudió sus trabajos.\\
El soporte público para las Tablas Rudolfinas se hizo casi
inmediatamente después de la publicación de Jacob Bartsch, un pupilo
de M\"{u}ller en Leipzig, posterior a (1630) sería un hermano-en la
Ley de Kepler-. En 1629 él publicó un volumen de Efemérides basado
sobre las tablas de Kepler, excepto que fueron calculadas para la
localidad de Strasbourg. En éstas él habló de las teorías de Kepler
en términos muy elogiosos, sólo que no son expandidas. Por tanto él
refiere a el lector al \textit{Epitome} para los principios teóricos
sobre los que las tablas están basadas.\\
Si ahora, avizoramos el período de 1630-1666, nos encontramos con
que el \textit{Epitome} Kepleriano atrajo modestamente la atención
cuando fue publicado por vez primera en 1618-1622, excepto que para
1630 fue evidente que estimuló un nuevo interés en el sistema
Copernicano.\\
Éste es fuertemente ilustrado por el hecho de que en 1631 dos
trabajos anticopernicanos fueron publicados por J.B.Morin en París y
Libert Froiddmont en Antwerp respectivamente. Ambos autores fueron
perturbados por el prestigioso aumento del Copernicanismo y ambos
hicieron del \textit{Epitome} uno de los dos blancos principales
para posteriores ataques- el otro es Landsberg`s
\textit{Commentationes in Motum Terrae}, publicado en 1630.
Froidmont da un breve bosquejo para las teorías físicas Keplerianas
excepto que escasamente mencionó las elipses. En 1633 Morin retornó
brevemente al mismo tema en el prefacio de su
\textit{Trigonometricae Canonicae}, donde él expresó el deseo de que
en un trabajo subsecuente tuviese la oportunidad de expandir sus
ideas sobre la teoría planetaria ''\textit{exmente Copernici et
Keppleri}''.\\
En 1632 las teorías Keplerianas fueron posteriormente soportadas por
el Astrónomo Alemán, Wilhelm Schickard.\\
En el año previo Gassendi, continúo una sugerencia de Kepler mismo,
y publicó una carta abierta para los astrónomos de Europa
preguntándoles a ellos por la observación del tránsito de Mercurio
hacia el Sol que tendría lugar el 7 de Noviembre de 1631. Schickard
fue uno de quienes cooperó. En un panfleto publicado en T\"{u}binge
(1632) él da un breve bosquejo para las ideas principales de Kepler,
incluyendo un estamento para las dos primeras leyes- la segunda es
únicamente mencionada, sin embargo, en la forma distancia-inversa
cualitativa. El refirió a los lectores al \textit{Epitome} y las
\textit{Rudolphine Tables} para detalles posteriores. Desde entonces
éste panfleto estuvo en forma de carta abierta para Gassendi quién
debería conocerlo, y, presumiblemente, para otros Astrónomos
Franceses.\\
En el mismo año, 1632, las dos primeras leyes de Kepler fueron
discutidas en el \textit{Directorium Generale} de Bonaventura
Cavalieri.\\
El \textit{Directorium} fue un libro de texto avanzado de astronomía
esférica, designado principalmente a la enseñanza del uso de los
logaritmos en computaciones astronómicas.\\
Pueden mencionarse otros dos usuarios más de las \textit{Rudolphine
Tables} por ésta época:\\
Adrian Vlacq en Holanda (1632), y Laurence Eichstadius en Stettin,
Alemania del Norte (1634).\\
En 1635 las teorías Keplerianas recibieron un estímulo posterior con
la publicación de una segunda Ed. para el \textit{Epitome} en
Frankfurt. Este trabajo es un volumen substancial de cerca de 100
pp; el hecho de que fuese reimpreso 5 años después de la muerte del
autor es una buena evidencia de un interés vivo en sus ideas durante
éste tiempo.\\
Otra referencia para las tablas de Kepler en 1630 puede se agregada.
Vincent Renieri, un monje Italiano y amigo de Galileo, publicó sus
\textit{Tabulae Medicaeae} en 1639, en la que da instrucciones
detalladas para el uso de seis conjuntos diferentes de tablas: las
de Kepler, Landsberg, Longomontanus, y otras tres antiguas: las
Pruténicas, Alphonsine y Ptolemaicas.\\
Las reglas fueron puramente prácticas; no había discusión sobre la
teoría y no se juraba sobre un uso relativo, exceptuando el hecho de
que aquí Kepler era colocado en primer lugar sugiriendo que
era el más importante de ellos. Claramente Renieri tuvo que estar
familiarizado por sí mismo con las dos primeras leyes.\\
En Bélgica el primer soporte prominente de Kepler fue el bien
conocido Astrónomo Gottefried Wendelin, quien en 1647 mencionó y
aceptó la primera Ley y en 1652 da un conteo detallado para la
tercera Ley. Él propuso, sin embargo, algunas modificaciones de la
primera Ley para la Luna, sugiriendo que su órbita fuese abolida en
lugar de la elíptica y que el foco primario no fuese exactamente el
centro de la Tierra, excepto que estuviese a 2,500 millas de
distancia de éste.\\
En 1650 encontramos una entusiasta discípula en María Cunitia- la
única mujer Astrónomo notable del siglo XVII. Una nativa de Silesia,
ella publicó su \textit{Urania Propitia} en Oels, 1650 (2 Ed.
Frankfurt, 1851). Este fue un conjunto simplificado de tablas
basadas sobre las Rudolphines y sobre las teorías físicas de Kepler,
con sus propias modificaciones para la segunda Ley. El volumen es
principalmente práctico en intención y ella da únicamente un breve
resumen para las ideas de Kepler, dejando en claro que ella las
aceptó de todo corazón.\\
Ahora, considerando que es de extremada importancia, hacer más
explícito la influencia Lutero-Platónica (como ya se ha dicho) en los
trabajos científicos del Matemático de los Estados de Estiria,
argüiremos que Kepler preocupado con la Astronomía Copernicana al
inicio de su carrera, desarrolló su sentido para la naturaleza
distintiva del Copernicanismo y las hipótesis astronómicas
únicamente en forma gradual. Un estado en tal desarrollo apareció en
un trabajo reluctantemente poco conocido, observándose que tuvo que
formar los puntos de vista de Kepler bajo la construcción de las
hipótesis astronómicas. En la puja de Tycho de Brahe, Kepler agregó
una composición como réplica contra Nicolaus Raymarus Ursus, quién
Tycho había acusado en 1596 de plagio. En 1597 Ursus afirmó en un
texto impreso que el sistema de Tycho no era de éste y atacó su
soberbia. Tycho respondió contra Ursus afirmando que éste era un
plagiador y llevó la acción hasta la corte. Desafortunadamente,
Kepler escribió en 1595 una carta laudatoria para Ursus, a quién él
en un tiempo llamó un matemático ejemplar de la época, y Ursus
incluyó esta cita en su carta para publicarla en la polémica contra
Tycho. Tycho incita a Kepler para que juntos escribieran una réplica
contra Ursus, una ''Apología Tychonis contra Ursus'', en la que ''la
naturaleza para las hipótesis astronómicas son explicadas y se
muestra que las afirmaciones de Ursus son incorrectas''. Este
tratado ''contra Ursus'' fue escrito en la década de 1600-1601
quedando incompleto debido al deceso de Tycho (24 de Octubre de
1601; Ursus había fallecido en Agosto de 1600) y no publicado hasta
la primera edición moderna de las obras de Kepler:
\textit{Johannes Kepler, Opera omnia, 8 Vols. ed. Cristian Frisch
(Frankfurt a. M: Heyder and Zimmer, 1858-1871), ''Contra Ursus''
Vol. I: 215-287.}\footnotemark\footnotetext{Cf. \textit{Kepler as
Historian of Science: Precursors of Copernican heliocentrism
according to Revolutionibus, I. 10.} Bruce Stansfield Eastwood.
Proceedings of the American Philosophical Society, Vol. 126. N. 5.
1982.}\\
Cuando escribiendo a Michael Maistlin (29 de Agosto de 1599), Kepler
habla de la afirmación absurda en Ursus de que la dirección tomada
por Tycho para su hipótesis fue pavimentada por Copérnico en el
\textit{De Revolutionibus} V, 35 (la explicación helio-estática para
las estaciones y retrogradaciones de Mercurio y Venus); y cuando
Ursus ridículamente encontró la hipótesis de Tycho en \textit{De
Revolutionibus} III, 25 (computando la aparente posición del Sol),
Kepler escribió ''asinus ille, nom jam Ursus''.\\
En breve, dice Kepler, Ursus no pudo haber leído a Ptolomeo ó tener
conocimiento de Copérnico. Por otro lado, después de haber pensado y
escrito sobre la materia (y después de la muerte de los
antagonistas), Kepler por sí mismo se muestra en un humor más
moderado. En una carta a David Fabricius (2 de Diciembre de 1602),
Kepler notó que deseaba estudiar a Proclo y Averroes sobre la
cuestión para las hipótesis antes de publicar ''contra Ursus'', y en
algún caso el deseaba observar el claro aire para el polvo de las
polémicas antes de imprimir el tratado. Entre esas dos cartas,
mientras escribía la defensa para Tycho, Kepler organizó sus
argumentos bajo cuatro capítulos anunciados:
\begin{enumerate}
\item Qué es una hipótesis Astronómica.\\
\item Una historia para la hipótesis Astronómica.\\
\item Una consideración para la hipótesis perdida de Apollonius.\\
\item La originalidad para la hipótesis de Tycho.
\end{enumerate}
En el interior del último de esos capítulos puede ser hallado un
fascinante repaso de Kepler para los posibles predecesores del
sistema planetario geo-heliocéntrico, entre quienes encontró a
Martianus Capella excepto que también a Macrobius, Pliny \index{Pliny}, Vitruvius \index{Vitruvius},
y hasta Platón. Este excursus sobre los orígenes de la idea
geo-heliocéntrica proporciona una comprensión del conocimiento de
Kepler para los precursores del Copernicanismo y su reconstrucción
para el significado de los textos involucrados.\\
Este es casi el lugar donde esos orígenes fueron mencionados. Kepler
deliberadamente los colocó a ellos como una afrenta para Ursus. No
únicamente falló Ursus en el conocimiento apropiado para los
escritos de Apollonius, que Ursus erróneamente afirmó como orígenes
de la hipótesis de Tycho.\\
El corazón para la referencia de Kepler en su repaso para el pasado
fue el \textit{De Revolutionibus} \index{De Revolutionibus}I, 10., en el que se lee, después
de ensayar la evidencia de los períodos orbitales para los planetas
en favor de un heliocentrismo razonable:
\begin{quote}\textit{''Consecuentemente Yo pienso que a nosotros no debería despistarnos
el argumento bien conocido de Martianus Capella, quien escribió la
\textit{Encyclopedia}, y otros ciertos escritores Latinos.\\
Ellos creían que Venus y Mercurio giraban alrededor del Sol que está
en el medio de ellos, y pensaban que esta es la razón por la que
ellos no divergen posteriormente de la curvatura para otras esferas
permitidas, por que ellos no giran alrededor de la Tierra, similar
al resto, excepto que tienen ápsides conversos (apsides conversae).
?`Qué otras personas, entonces, comprenderían, que el centro para
las órbitas está en la región del
Sol?}.\footnotemark\footnotetext{Cf. Ibid. pp. 368. n. 9.}\end{quote}

Los trabajos de Platón \index{Platón} sobre los planetas inertes también son
obscuros, las observaciones de Kepler, son pistas que muestran su
significado real como Macrobius posteriormente explicará. El caso es
hecho por Platón por el empleo de sus propias palabras; Martianus, y
especialmente Macrobius. Kepler inicia con una cita del
\textit{Timaeus} 38 C-D:\\
''Cuando Dios tuvo que hacer los diferentes cuerpos [trotamundos],
el conjunto de sus revoluciones con el circuito que cada uno
requeriría, seis [cuerpos] en seis [órbitas], la Luna en el primero
a continuación la Tierra, el Sol en el segundo bajo la Tierra,
entonces la Estrella Nocturna [Venus] y el uno sagrado llamado
Mercurio en tales [círculos] como por la razón de sus velocidades
iguales para la órbita del Sol, aunque hayan recibido una potencia
contraria (contrariam vim) para este, por tanto el Sol, Mercurio y
Venus incluido (comprehendunt) están incluidos uno en cada
otro.''\footnotemark\footnotetext{Cf. Ibid. pp. 381-82. n.
74-75.}\\
Es notable que Kepler empleara únicamente el pasaje del
\textit{Timaeus} y no el pasaje corroborativo en la República. Él
pudo haber elegido emplear únicamente el trabajo conocido durante la
Edad Media Latina. Él lo elige por que el \textit{Timaeus} contiene
contiene el pensamiento cosmológico, y la República no lo hace,
observándose que es el origen apropiado para un tema cosmológico.\\
Ahora, es conveniente subrayar, que para Copérnico Platón
simplemente en el \textit{Timaeus} coloca a Mercurio y Venus bajo el
Sol.\footnotemark\footnotetext{\textit{De Revolutionibus} (1543),
f.7v,19-22. transl. Rosen, \textit{On the Revolutions}, pp. 18,
35-36. Cf. Eastwood. Op. cit. n. 78.}\\
En el establecimiento del pedigree Platónico, Kepler usa la
narración de Macrobius como pieza central para la reconstrucción del
pensamiento antiguo.\\
Tomando una sucesión de citas del \textit{Commentary on the Dream of
Scipio}, I, 19 de Macrobius, Kepler hace las siguientes anotaciones:
\begin{enumerate}
\item Ambas tradiciones la supra-solar y la infra-solar
concernientes a la posición de Mercurio y Venus son correctas, y
Macrobius conocía esto.\\
La causa para la inversión intermitente de Mercurio y Venus no es
buscada, y Macrobius dice que Platón tuvo una explicación.\\
\item Platón razonó que el orden para los planetas corresponde a su
período orbital.\\
\item En la presentación de las órbitas de Venus, Mercurio y el Sol,
Macrobius habla de los tres teniendo un simple coelum en el que
ellos orbitan en un período de, más o menos, un año. Esto Kepler lo
encuentra significativo, el coelum tendrá una física distinta y
no justamente un fenómeno de igualdad- como sí los tres estuvieran
empatados y simplemente no aparecieran juntos.\\
\item Los planetas Mercurio y Venus están mucho más propiamente
asociados con el Sol que con la Luna, los movimientos son más
cerrados en el Sol que en la Luna.\\
\item El sentido completo de Platón interpretado en la dirección de
Macrobio puede ser auxiliado por la mejoría de una simple palabra,
dice Kepler.
\end{enumerate}
En muchas direcciones el \textit{Commentary} de Macrobius constituye
el testigo clave para Kepler en la traza de la conexión Platónica.\\
Para todos los orígenes Keplerianos Macrobius \index{Macrobius} es el único autor que
reconoce el Platonismo y el Heliocentrismo.\\
Dos términos usados en el pasaje de Kepler en ''contra Ursus'' y
requiriendo delimitaciones son \textit{circulus} y \textit{vertex}.
Macrobius emplea circulus como órbita, no como epiciclo. Para
suponer la interpretación de los epiciclos se requiere un número de
suposiciones inaceptables. Esto requiere primero que todo que
Macrobius leyó el \textit{Commentary on Platon's Timaeus} de
Calcidius \index{Calcidius} y considerara involucrar la explicación de los epiciclos
hallado en el \textit{Commentary}.\\
Aún cuando el \textit{Commentary} de
Calcidius\footnotemark\footnotetext{En la versión de Calcidius para
el Timaeus 36 D (Ed. Waszink, pp. 28, 21-24) el movimiento relativo
para cada planeta no es dado; uno encuentra únicamente que tres
movimientos para la misma velocidad y las cuatro velocidades son
proporcionales. Esas traslaciones por Calcidius representan
adecuadamente a Platón, así que se necesitan observar los
comentarios sobre Platón para más definiciones. El comentario de
Calcidius sobre el Timaeus \index{Timaeus} ofrece una clarificación no explícita
para el orden exacto del Sol-Mercurio-Venus elaborada en 38 D; y
nosotros aprendemos el orden detallado únicamente de su comentario
sobre 36 D, donde las proporciones para el movimiento planetario son
dadas. Calcidius emplea la sucesión 1,2,3,4,9,8,27 como múltiplos
para la distancia unidad Tierra-Luna en la localización para las
posiciones relativas de los planetas. Macrobius sigue un trabajo
perdido de Porphyry, no a Calcidius, en el empleo de la sucesión
1,2,3,4,9,8,27 como múltiplos sucesivos cada uno para la distancia
planetaria inmediatamente precedente, así Marte, V.g., es 9 veces la
distancia entre la Tierra y Mercurio. únicamente con esta
elaboración posterior de Platón por los comentaristas Platónicos,
tales como Calcidius y Porphyry, el Timaeus da una definición
adecuada para el orden adecuado de los planetas, y en todos los
comentadores se observa ignoran la dificultad en la afirmación de
que el Sol con Venus y Mercurio tienen períodos orbitales similares
mientras se observan tres marcas diferentes en los intervalos
''armónicos''. Cf, Ibid. n. 90.} emplee círculos comúnmente para
órbitas ó órbitas aparentemente en lugar de epiciclos. Sin embargo,
Calcidius describe distintivamente el orden Platónico para los
planetas como ellos son dados por Platón, que se comportan de forma
diferente en Macrobius el orden Platónico para la
Luna-Mercurio-Venus, y así. La interpretación para los círculos de
Macrobius como epiciclos requiere la suposición de que su audiencia,
ostensiblemente su Joven hijo, Eustachius, y más generalmente una
audiencia no-sofisticada (astronómicamente), inmediatamente
conocieran el futuro del círculos (órbita) para el Sol del círculos
(epiciclo) para el planeta, aunque Macrobius nunca mencionara los
epiciclos en algún lugar del \textit{Commentary}!.\\
?`Cómo hace Kepler para encontrar una trayectoria heliocéntrica en
el pasaje de Macrobius? En primera instancia se observa que él
conocía en avanzada lo que deseaba ó esperaba hallar. Evidencia
colateral ayuda al marco de sus oportunidades. Por otra parte,
ninguna de las referencia tempranas en el siglo XVI para el sistema
de Martianus hace alguna referencia a Macrobius, ó las referencias
del siglo XVII tales como Riccioli, Argoli, Sherburne y Cellarius.
De otro lado, las adiciones de Macrobius son una referencia para
Martianus Capella en tiempo de Kepler. Para las 33 ediciones
separadas antes de 1600 para el \textit{Commentary of the Dream of
Scipio} de Macrobius se tiene una, tal como la de Leiden, 1597,
añade una nota para el estamento de Macrobius concerniente a las
órbitas (circuli) para el Sol-Mercurio-Venus es equivalente a
Copérnico para el conocimiento de Martianus Capella. Es conveniente
agregar que Kepler hizo esta anotación, excepto que no hay razón
para dudar de esto. Sobre la base de las citas de Kepler para
Macrobius, el texto de Leiden 1597 se observa que no pudo ser uno de
los orígenes de Kepler, al menos cuando se asienta bajo la copia del
texto de Macrobius.\\
El último origen Kepleriano para el esclarecimiento sobre Platón es
el escrito identificado como el \textit{Vetus commentarius} de
Bede.\\
Las citas hechas dejan en claro que éste \textit{Vetus commentarius}
emplea el comentario de Calcidius sobre el Timaeus de Platón,
excepto que Kepler no lo reconoce e identifica el trabajo en tal
dirección.

\end{document}